\input amstex
\input amsppt.sty
\magnification=\magstep1
%\hsize=32truecc
%\vsize=22.2truecm
\hsize=33.5truecc
\vsize=23truecm
\baselineskip=16truept
\NoBlackBoxes
\TagsOnRight \pageno=1 \nologo
\def\Z{\Bbb Z}
\def\N{\Bbb N}

\def\l{\left}
\def\r{\right}
\def\bg{\bigg}
\def\({\bg(}
\def\[{\bg\lfloor}
\def\){\bg)}
\def\]{\bg\rfloor}
\def\t{\text}
\def\f{\frac}

\def\p{\ (\roman{mod}\ p)}

\def\sm{\setminus}

\def\bi{\binom}
\def\eq{\equiv}

\def\ls{\leqslant}
\def\gs{\geqslant}
\def\mo{\roman{mod}}
\def\ord{\roman{ord}}

\def\al{\alpha}
\def\da{\delta}

\def\Proof{\noindent{\it Proof}}

\def\Remark{\medskip\noindent{\it  Remark}}

\hbox {Preprint, {\tt arXiv:1610.03384}}
\bigskip
\topmatter
\title Supercongruences involving Lucas sequences\endtitle
\author Zhi-Wei Sun\endauthor
\leftheadtext{Zhi-Wei Sun}
\affil Department of Mathematics, Nanjing University\\
 Nanjing 210093, People's Republic of China
  \\  zwsun\@nju.edu.cn
  \\ {\tt http://maths.nju.edu.cn/$\sim$zwsun}
\endaffil
\abstract For $A,B\in\Z$, the Lucas sequence $u_n(A,B)\ (n=0,1,2,\ldots)$ are defined by
$u_0(A,B)=0$, $u_1(A,B)=1$, and $u_{n+1}(A,B)=Au_n(A,B)-Bu_{n-1}(A,B)\ (n=1,2,3,\ldots).$
For any odd prime $p$ and positive integer $n$, we establish the new result
$$\f{u_{pn}(A,B)-(\f{A^2-4B}p)u_n(A,B)}{pn}\in\Z_p,$$
where $(\f{\cdot}p)$ is the Legendre symbol and $\Z_p$ is the ring of $p$-adic integers.

Let $p$ be an odd prime and let $n$ be a positive integer. For any integer $m\not\equiv0\pmod p$, we show
that $$\f1{pn}\bigg(\sum_{k=0}^{pn-1}\f{\bi{2k}k}{m^k}-\l(\f{\Delta}p\r)\sum_{r=0}^{n-1}\f{\bi{2r}r}{m^r}\bigg)\in\Z_p$$
and furthermore
$$\f1n\bigg(\sum_{k=0}^{pn-1}\f{\bi{2k}k}{m^k}-\l(\f{\Delta}p\r)\sum_{r=0}^{n-1}\f{\bi{2r}r}{m^r}\bigg)
\equiv\f{\bi{2n-1}{n-1}}{m^{n-1}}u_{p-(\f{\Delta}p)}(m-2,1)\pmod{p^2},$$
where $\Delta=m(m-4)$.
We also pose some conjectures for further research.
\endabstract
\thanks 2020 {\it Mathematics Subject Classification}. \,Primary 11A07, 11B65;
Secondary  05A10, 11B39, 11B75.
\newline\indent Supported by the National Natural Science
Foundation (grant 11971222) of China.
\endthanks

\endtopmatter
\document

\heading{1. Introduction}\endheading

Let $p>3$ be a prime. In 2006 H. Pan and the author [PS] deduced from a sophisticated combinatorial identity the congruence
$$\sum_{k=0}^{p-1}\bi{2k}{k+d}\eq\l(\f {p-d}3\r)\pmod p\ \quad \t{for}\ d=0,\ldots,p-1,$$
where $(-)$ denotes the Legendre symbol.
In 2011 the author and R. Tauraso [ST11] showed further that
$$\sum_{k=0}^{p-1}\bi{2k}{k+d}\eq\l(\f {p-d}3\r)\pmod{p^2}\quad\ \t{for}\ d=0,1.\tag1.1$$
Recently, J.-C. Liu [L16] proved the following extension with $n\in\Z^+=\{1,2,3,\ldots\}$ conjectured by M. Apagodu and D. Zeilberger [AZ]:
$$\sum_{k=0}^{pn-1}\bi{2k}k\eq\l(\f p3\r)\sum_{k=0}^{n-1}\bi{2k}k\pmod{p^2}\tag1.2$$
and
$$\aligned\sum_{k=0}^{pn-1}C_k\eq\cases\sum_{r=0}^{n-1}C_r\pmod{p^2}&\t{if}\ p\eq1\pmod3,
\\-\sum_{r=0}^{n-1}(3r+2)C_r\pmod{p^2}&\t{if}\ p\eq2\pmod3.
\endcases\endaligned\tag1.3$$
where $C_k$ denotes the Catalan number $\bi{2k}k-\bi{2k}{k+1}=\bi{2k}k/(k+1)$.
Note that this result in the case $n=1$ yields the supercongruence (1.1).

For given integers $A$ and $B$, the Lucas sequence $u_n=u_n(A,B)\ (n\in\N=\{0,1,2,\ldots\})$ is given by
$$u_0=0,\ u_1=1,\ \quad\t{and}\ u_{n+1}=Au_{n}-Bu_{n-1}\ \t{for}\ n=1,2,3,\ldots.$$
It is well known that $p\mid u_{p-(\f{A^2-4B}p)}$ for any odd prime $p$ not dividing $B$ (see, e.g., [S10, Lemma 2.3]).
In 2010 the author [S10] showed that for any nonzero integer $m$ and odd prime $p$ not dividing $m$ we have
$$\sum_{k=0}^{p-1}\f{\bi{2k}k}{m^k}\eq\l(\f{m(m-4)}p\r)+u_{p-(\f{m(m-4)}p)}(m-2,1)\pmod{p^2}.\tag1.4$$

In this paper we obtain the following general result which is a common extension of (1.1)-(1.4).

\proclaim{Theorem 1.1} Let $m\in\Z\sm\{0\}$, $n\in\Z^+$ and $\Delta=m(m-4)$. For any odd prime $p$ not dividing $m$, we have
$$\f1n\(\sum_{k=0}^{pn-1}\f{\bi{2k}k}{m^k}-\l(\f{\Delta}p\r)\sum_{r=0}^{n-1}\f{\bi{2r}r}{m^r}\)
\eq\f{\bi{2n-1}{n-1}}{m^{n-1}}u_{p-(\f{\Delta}p)}(m-2,1)\pmod{p^2}\tag1.5$$
and
$$\aligned&\f1n\(\sum_{k=0}^{pn-1}\f{\bi{2k}{k+1}}{m^k}-\l(\f{\Delta}p\r)\sum_{r=0}^{n-1}\f{\bi{2r}{r+1}}{m^r}
+\l(\f m2-\f{\bi{2n-1}{n-1}}{m^{n-1}}\r)\l(1-\l(\f{\Delta}p\r)\r)\)
\\\eq&\f{\bi{2n-1}{n-1}}{m^{n-1}}\l(1-m^{p-1}+\f{m-2}{2}u_{p-(\f{\Delta}p)}(m-2,1)\r)\pmod{p^2},
\endaligned\tag1.6$$
hence
$$\aligned&\f1n\(\sum_{k=0}^{pn-1}\f{C_k}{m^k}-\l(\f{\Delta}p\r)\sum_{r=0}^{n-1}\f{C_r}{m^r}
+\l(\f{\bi{2n-1}{n-1}}{m^{n-1}}-\f m2\r)\l(1-\l(\f{\Delta}p\r)\r)\)
\\\eq&\f{\bi{2n-1}{n-1}}{m^{n-1}}\l(m^{p-1}-1+\f{4-m}{2}u_{p-(\f{\Delta}p)}(m-2,1)\r)\pmod{p^2}.
\endaligned\tag1.7$$
\endproclaim

\proclaim{Corollary 1.1} Let $p$ be an odd prime and let $n\in\Z^+$. Then
$$\align\f1n\(\sum_{k=0}^{pn-1}\bi{2k}k-\l(\f{p}3\r)\sum_{r=0}^{n-1}\bi{2r}r\)\eq&0\pmod{p^2},\tag1.8
\\\f1n\(\sum_{k=0}^{pn-1}\f{\bi{2k}k}{2^k}-\l(\f{-1}p\r)\sum_{r=0}^{n-1}\f{\bi{2r}r}{2^r}\)\eq&0\pmod{p^2},\tag1.9
\endalign$$
and
$$\f1{n}\(\sum_{k=0}^{pn-1}C_k-\l(\f p3\r)\sum_{r=0}^{n-1}C_r+\f{1-(\f p3)}2\l(\bi{2n}n-1\r)\)\eq0\pmod{p^2}.\tag1.10$$
When $p>3$, we also have
$$\f1n\(\sum_{k=0}^{pn-1}\f{\bi{2k}k}{3^k}-\l(\f p3\r)\sum_{k=0}^{n-1}\f{\bi{2r}r}{3^r}\)\eq0\pmod{p^2}.\tag1.11$$
\endproclaim
\Proof. By induction,  $u_{2k}(0,1)=0$ and $u_{3k}(\pm1,1)=0$ for all $k\in\N$.
Applying (1.5) with $m=1,2,3$  we obtain (1.8), (1.9) and (1.11).
(1.7) with $m=1$ yields (1.10).
 \qed

\Remark\ 1.1. Our (1.11) implies (1.3) since $\sum_{r=0}^{n-1}(3r+1)C_r=\bi{2n}n-1$ for all $n\in\Z^+$.
(1.9) and (1.10) in the case $n=1$ were first proved by the author [S11b].
\medskip

For given integers $A$ and $B$, the sequence $v_n=v_n(A,B)\ (n=0,1,2,\ldots)$ defined by
$$v_0=2,\ v_1=A,\ \t{and}\ v_{n+1}=Av_n-Bv_{n-1}\ (n=1,2,3,\ldots)$$
is called the the companion sequence of the Lucas sequence $u_n=u_n(A,B)\ (n\in\N)$.
By induction,
$$v_{n}(A,B)=2u_{n+1}(A,B)-Au_n(A,B)\quad \ \t{for all}\ n\in\N.$$

To prove Theorem 1.1, we need the following auxiliary result on general Lucas sequences which has its own interest.

\proclaim{Theorem 1.2} Let $A,B\in\Z$ and $\Delta=A^2-4B$. Let $p$ be an odd prime and let $n\in\Z^+$. Then
$$\f{u_{pn}(A,B)-(\f{\Delta}p)u_n(A,B)}{pn}\in\Z_p\ \ \t{and}\ \ \f{v_{pn}(A,B)-v_n(A,B)}{pn}\in\Z_p,\tag1.12$$
where $\Z_p$ denotes the ring of $p$-adic integers. Moreover, if $p\nmid B\Delta$ then
$$\aligned\f{u_{pn}(A,B)-(\f{\Delta}p)u_n(A,B)}{pn}
\eq&\f{u_n(A,B)}2\l(\f{\Delta}p\r)\f{B^{p-1}-1}p
\\&+\f{v_n(A,B)}{2B^{(1-(\f{\Delta}p))/2}}\cdot\f{u_{p-(\f{\Delta}p)}(A,B)}p
\pmod{p}
\endaligned\tag1.13$$
and
$$\aligned\f{v_{pn}(A,B)-v_n(A,B)}{pn}
\eq&\f{v_n(A,B)}2\cdot\f{B^{p-1}-1}p
\\&+\f{\Delta u_n(A,B)}{2B^{(1-(\f{\Delta}p))/2}}\l(\f{\Delta}p\r)\f{u_{p-(\f{\Delta}p)}(A,B)}p
\pmod{p}.\endaligned\tag1.14$$
\endproclaim
\Remark\ 1.2. (1.12) in the case $n=1$ is well known, see, e.g., [S10, Lemma 2.3]. For the prime $p=2$, (1.12) also holds
if we adopt the Kronecker symbol
$$\l(\f{\Delta}2\r)=\cases1&\t{if}\ \Delta\eq1\pmod8,
\\-1&\t{if}\ \Delta\eq5\pmod 8,
\\0&\t{if}\ \Delta\eq0\pmod2.
\endcases$$

Motivated by Theorems 1.1 and 1.2, we give another theorem on supercongruences.

\proclaim{Theorem 1.3}
For any prime $p>5$ and $n\in\Z^+$, we have
$$\f{g_{pn}(-1)-g_n(-1)}{n^2}\eq0\pmod {p^3},\tag1.15$$
where
$$g_n(x):=\sum_{k=0}^n\bi nk^2\bi{2k}kx^k.\tag1.16$$
\endproclaim
\Remark\ 1.3. The polynomials $g_n(x)$ $(n=0,1,2,\ldots)$ were introduced by the author [S16] in which the author proved for any prime $p>5$ that
$$\sum_{k=1}^{p-1}\f{g_k(-1)}{k}\eq0\pmod{p^2}\ \ \t{and}\ \ \sum_{k=1}^{p-1}\f{g_k(-1)}{k^2}\eq0\pmod{p};$$
see also V.J. Guo, G.-S. Mao and H. Pan [GMP] for some congruences involving the polynomials $g_n(x)\ (n=0,1,2,\ldots)$.
\medskip

We will show Theorems 1.2, 1.1 and 1.3 in Sections 2, 3 and 4 respectively. In Section 5 we pose some conjectures for further research.

\heading{2. Proof of Theorem 1.2}\endheading

Let $A,B\in\Z$, and let
 $$\al=\f{A+\sqrt{\Delta}}2\ \ \t{and}\ \ \beta=\f{A-\sqrt{\Delta}}2\tag2.1$$
 be the two roots of the quadratic equation $x^2-Ax+B=0$, where $\Delta=A^2-4B$.
 It is well known that
 $$(\al-\beta)u_n(A,B)=\al^n-\beta^n\ \ \t{and}\ \ v_n(A,B)=\al^n+\beta^n\quad\t{for all}\ n\in\N.\tag2.2$$
 When $\Delta=0$, by induction we have $u_n(A,B)=n(A/2)^{n-1}$ for all $n\in\Z^+$.

\proclaim{Lemma 2.1} Let $A,B\in\Z$ and $n\in\Z^+$. Then
$$u_n(A,B)=\sum_{k=0}^{\lfloor(n-1)/2\rfloor}\bi{n-1-k}kA^{n-1-2k}(-B)^k\tag2.3$$
and
$$v_n(A,B)=\sum_{k=0}^{\lfloor n/2\rfloor}\f n{n-k}\bi{n-k}kA^{n-2k}(-B)^k.\tag2.4$$
\endproclaim
\Remark\ 2.1. Lemma 2.1 is a well known result (see, e.g., [G, (1.60)]) and the two identities (2.3) and (2.4) can be easily proved by induction on $n$.

\proclaim{Lemma 3.2 {\rm [S11a, Lemma 2.2]}} Let $A,B\in\Z$ and let $d\in\Z^+$ be an odd divisor of $\Delta=A^2-4B$.
Then, for any $n\in\Z^+$, we have
$$\f{u_n(A,B)}n\eq\l(\f A2\r)^{n-1}+\cases(A/2)^{n-3}\Delta/3\ (\mo\ d)&\t{if}\ 3\mid d\ \t{and}\ 3\mid n,
\\0\ (\mo\ d)&\t{otherwise}.\endcases\tag2.5$$
\endproclaim

\proclaim{Lemma 2.3 {\rm [S12b, Lemma 2.2]}}  Let $A,B\in\Z$ and $\Delta=A^2-4B$.
Suppose that $p$ is an odd prime with $p\nmid B\Delta$. Then we have the congruence
$$\l(\f{A\pm\sqrt{\Delta}}2\r)^{p-(\f{\Delta}p)}\eq B^{(1-(\f{\Delta}p))/2}\ \ (\mo\ p)\tag2.6$$
in the ring of algebraic integers.
\endproclaim

\proclaim{Lemma 2.4} Let $A,B\in\Z$ and $\Delta=A^2-4B$. For any odd prime $p$ not dividing $B\Delta$, we have
$$v_{p-(\f{\Delta}p)}(A,B)\eq B^{(1-(\f{\Delta}p))/2}(B^{p-1}+1)\pmod{p^2}.\tag2.7$$
\endproclaim
\Proof. By the proof of [S13, Lemma 3.2],
$$v_{p-(\f{\Delta}p)}(A,B)\eq2\l(\f Bp\r)B^{(p-(\f{\Delta}p))/2}\pmod{p^2}.$$
Thus
$$\align &v_{p-(\f{\Delta}p)}(A,B)-B^{(1-(\f{\Delta}p))/2}(B^{p-1}+1)
\\\eq&2\l(\f Bp\r)B^{(p-(\f{\Delta}p))/2}-B^{(1-(\f{\Delta}p))/2}(B^{p-1}+1)
\\=&B^{(1-(\f{\Delta}p))/2}\l(2\l(\l(\f{B}p\r)B^{(p-1)/2}-1\r)-(B^{p-1}-1)\r)
\\\eq&0\pmod{p^2}
\endalign$$
since
$$\align B^{p-1}-1=&\l(\l(\f Bp\r)B^{(p-1)/2}+1\r)\l(\l(\f Bp\r)B^{(p-1)/2}-1\r)
\\\eq&2\l(\l(\f Bp\r)B^{(p-1)/2}-1\r)\pmod{p^2}.
\endalign$$
This concludes the proof. \qed

\proclaim{Lemma 2.5} Let $p$ be a prime and let $n\in\Z^+$. For any $p$-adic integer $a\not\eq0\pmod p$ and positive integer $n$, we have
$$\f{a^{(p-1)n}-1}{pn}\in\Z_p\tag2.8$$
and moreover
$$\f{a^{(p-1)n}-1}{pn}\eq n^{\da_{p,2}}\f{a^{p-1}-1}p\pmod p.\tag2.9$$
\endproclaim
\Proof. Let $r$ be the unique integer in $\{1,\ldots,p-1\}$ with $a\eq r\pmod p$. Then $a^{p-1}\eq r^{p-1}\eq1\pmod p$ by Fermat's little theorem.
Write $a^{p-1}=1+pt$ with $t\in\Z_p$. Observe that
$$\align\f{a^{(p-1)n}-1}{pn}=&\f{(1+pt)^n-1}{pn}=\f1{pn}\sum_{k=1}^n\bi nk(pt)^k=\sum_{k=1}^n\bi{n-1}{k-1}\f{p^{k-1}}kt^k
\\\eq&t+(n-1)\f p2t^2\eq t+(n-1)\f p2t\eq n^{\da_{p,2}}t\pmod p
\endalign$$
since $p^{k-2}/k\in\Z_p$ for all $k=3,4,\ldots$. This concludes the proof. \qed

\medskip
\noindent{\it Proof of Theorem 1.2}. For the sake of brevity, we just write $u_k=u_k(A,B)$ and $v_k=v_k(A,B)$ for all $k\in\N$.
Let $\al$ and $\beta$ be the algebraic integers defined by (2.1).

 If $p\mid \Delta$, then by Lemma 2.2 we have
$$\f{u_{pn}-(\f{\Delta}p)u_n}{pn}=\f{u_{pn}}{pn}\eq\l(\f{A}2\r)^{pn-1}+\da_{p,3}\l(\f{A}2\r)^{3n-3}\f{\Delta}3\pmod p.$$
By (2.2),
$$v_n=\l(\f{A+\sqrt{\Delta}}2\r)^n+\l(\f{A-\sqrt{\Delta}}2\r)^n=\f1{2^{n-1}}\sum_{k=0}^{\lfloor n/2\rfloor}\bi n{2k}A^{n-2k}\Delta^k.$$
When $p\mid \Delta$, we have
$$\f{v_n-A^n/2^{n-1}}{pn}=\f{\Delta}p\sum_{0<k\ls\lfloor n/2\rfloor}\bi{n-1}{2k-1}A^{n-2k}\f{\Delta^{k-1}}{k2^n}\in\Z_p$$
and similarly
$$\f{v_{pn}-A^{pn}/2^{pn-1}}{pn}\in\Z_p,$$
hence $(v_{pn}-v_n)/(pn)\in\Z_p$ since
$$\f{A^{pn}/2^{pn-1}-A^n/2^{n-1}}{pn}=\f{A^n}{2^{n-1}pn}\l(\l(\f A2\r)^{(p-1)n}-1\r)\in\Z_p$$
by Lemma 2.5.

 Below we assume that $p\nmid \Delta$. Note that
 $$\align u_{pn}=&\f{\alpha^{pn}-\beta^{pn}}{\al-\beta}=\f{\al^n-\beta^n}{\al-\beta}\cdot\f{(\al^n)^p-(\beta^n)^p}{\al^n-\beta^n}
 \\=&u_nu_p(\al^n+\beta^n,\al^n\beta^n)=u_nu_p(v_n,B^n)
 \endalign$$
 and
 $$v_{pn}=(\al^n)^p+(\beta^n)^p=v_p(\al^n+\beta^n,\al^n\beta^n)=v_p(v_n,B^n).$$
 By Lemma 2.1,
 $$u_p(v_n,B^n)=\sum_{k=0}^{(p-1)/2}\bi {p-1-k}kv_n^{p-1-2k}(-B^n)^k$$
 and
 $$v_p(v_n,B^n)=\sum_{k=0}^{(p-1)/2}\f p{p-k}\bi {p-k}kv_n^{p-2k}(-B^n)^k.$$

 Now suppose that $p\mid B$. Then $(-B^n)^k/(pn)\in\Z_p$ for all $k\in\Z^+$ since $p^{n-1}/n\in\Z_p$.
 Note also that $(\f{\Delta}p)=(\f{A^2}p)=1$. Thus
 $$\f{u_{pn}-(\f{\Delta}p)u_n}{pn}-\f{u_n}{pn}(v_n^{p-1}-1)\in\Z_p\ \ \t{and}\ \ \f{v_{pn}-v_n^p}{pn}\in\Z_p.$$
 In view of (2.4),
 $$v_n-A^n=-Bn\sum_{0<k\ls\lfloor n/2\rfloor}\bi{n-k-1}{k-1}\f{(-B)^{k-1}}{k}A^{n-2k}.$$
 Since $p\mid B$ and $p^{k-1}/k\in\Z_p$ for all $k\in\Z^+$, we see that $v_n=A^n+pnt$ for some $t\in\Z_p$.
 By Lemma 2.5, $(A^{(p-1)n}-1)/(pn)\in\Z_p$. Therefore $(v_n^{p-1}-1)/(pn)\in\Z_p$
 and hence $(u_{pn}-(\f{\Delta}p)u_n)/(pn)\in\Z_p$. Note also that
 $$\f{v_{pn}-v_n}{pn}=\f{v_{pn}-v_n^p}{pn}+v_n\f{v_n^{p-1}-1}{pn}\in\Z_p.$$

Below we suppose that $p\nmid B\Delta$.

{\it Case} 1. $(\f{\Delta}p)=1$.

In this case, by Lemma 2.3 we have $\alpha^{p-1}\eq\beta^{p-1}\eq1\pmod p$ in the ring of algebraic integers.
Similar to Lemma 2.5, we have
$$\f{\al^{(p-1)n}-1}{pn}\eq\f{\al^{p-1}-1}p\pmod p\ \ \t{and}\ \ \f{\beta^{(p-1)n}-1}{pn}\eq\f{\beta^{p-1}-1}p\pmod p.$$
 Therefore
$$\align\f{u_{pn}-(\f{\Delta}p)u_n}{pn}=&\f{\al^{pn}-\beta^{pn}-(\al^n-\beta^n)}{pn(\al-\beta)}
\\=&\f{\al-\beta}{\Delta}\l(\al^n\f{\al^{(p-1)n}-1}{pn}-\beta^n\f{\beta^{(p-1)n}-1}{pn}\r)
\\\eq&\f1{\al-\beta}\cdot\f{\al^n(\al^{p-1}-1)-\beta^n(\beta^{p-1}-1)}{p}
\\=&\f{(\al^n-\beta^n)(\al^{p-1}+\beta^{p-1}-2)+(\al^n+\beta^n)(\al^{p-1}-\beta^{p-1})}{2p(\al-\beta)}
\\=&\f{u_n}2\cdot\f{v_{p-1}-2}p+\f{v_n}2\cdot\f{u_{p-1}}p\pmod p\endalign$$
and
$$\align\f{v_{pn}-v_n}{pn}=&\f{\al^{pn}+\beta^{pn}-(\al^n+\beta^n)}{pn}
\\=&\al^n\f{\al^{(p-1)n}-1}{pn}+\beta^n\f{\beta^{(p-1)n}-1}{pn}
\\\eq&\f{\al^n(\al^{p-1}-1)+\beta^n(\beta^{p-1}-1)}p
\\=&\f{(\al^n+\beta^n)(\al^{p-1}+\beta^{p-1}-2)+(\al^n-\beta^n)(\al^{p-1}-\beta^{p-1})}{2p}
\\=&\f{v_n}2\cdot\f{v_{p-1}-2}p+\f{\Delta u_n}2\cdot\f{u_{p-1}}p\pmod p.
\endalign$$

{\it Case} 2. $(\f{\Delta}p)=-1$.

In this case, by Lemma 2.3 we have $\alpha^{p+1}\eq\beta^{p+1}\eq B\pmod p$ in the ring of algebraic integers.
Thus
$$\align\f{\al^{(p+1)n}-B^n}{pn}=&\f1{pn}\sum_{k=1}^n\bi nk(\al^{p+1}-B)^kB^{n-k}
\\=&\f{\al^{p+1}-B}p\sum_{k=1}^n\bi{n-1}{k-1}\f{(\al^{p+1}-B)^{k-1}}{k}B^{n-k}
\\\eq& \f{\al^{p+1}-B}pB^{n-1}\pmod p.\endalign$$
Similarly,
$$\f{\beta^{(p+1)n}-B^n}{pn}\eq \f{\beta^{p+1}-B}pB^{n-1}\pmod p.$$ Therefore
$$\align\f{u_{pn}-(\f{\Delta}p)u_n}{pn}=&\f{\al^{pn}-\beta^{pn}+\al^n-\beta^n}{pn(\al-\beta)}
\\=&\f{\al-\beta}{pnB^n\Delta}\l((\al\beta)^n\al^{pn}-(\al\beta)^n\beta^{pn}+B^n(\al^n-\beta^n)\r)
\\=&\f{\al-\beta}{B^n\Delta}\l(\beta^n\f{\al^{(p+1)n}-B^n}{pn}+\al^n\f{B^n-\beta^{(p-1)n}}{pn}\r)
\\\eq&\f1{(\al-\beta)B}\l(\beta^n\f{\al^{p+1}-B}p+\al^n\f{B-\beta^{p+1}}p\r)
\\=&\f{v_n}{2B}\cdot\f{u_{p+1}}p-\f{u_n}{2B}\cdot\f{v_{p+1}-2B}p\pmod p
\endalign$$
and
$$\align \f{v_{pn}-v_n}{pn}=&\f{(\al\beta)^n}{pnB^n}(\al^{pn}+\beta^{pn})-\f{\al^n+\beta^n}{pn}
\\=&\f1{B^n}\l(\beta^n\f{\al^{(p+1)n}-B^n}{pn}+\al^n\f{\beta^{(p+1)n}-B^n}{pn}\r)
\\\eq&\f1{B}\l(\beta^n\f{\al^{p+1}-B}p+\al^n\f{\beta^{p+1}-B}p\r)
\\=&\f{(\al^n+\beta^n)(\al^{p+1}+\beta^{p+1}-2B)-(\al^n-\beta^n)(\al^{p+1}-\beta^{p+1})}{2Bp}
\\=&\f{v_n}{2B}\cdot\f{v_{p+1}-2B}p-\f{\Delta u_n}{2B}\cdot\f{u_{p+1}}p\pmod p.
\endalign$$

Whether $(\f{\Delta}p)$ is $1$ or $-1$, we always have
$$\align\f{u_{pn}-(\f{\Delta}p)u_n}{pn}\eq&\l(\f{\Delta}p\r)\f{u_n}{2B^{(1-(\f{\Delta}p))/2}}\cdot\f{v_{p-(\f{\Delta}p)}-2B^{(1-(\f{\Delta}p))/2}}p
\\&+\f{v_n}{2B^{(1-(\f{\Delta}p))/2}}\cdot\f{u_{p-(\f{\Delta}p)}}p\pmod p
\endalign$$
and
$$\align\f{v_{pn}-v_n}{pn}\eq&\f{v_n}{2B^{(1-(\f{\Delta}p))/2}}\cdot\f{v_{p-(\f{\Delta}p)}-2B^{(1-(\f{\Delta}p))/2}}p
\\&+\l(\f{\Delta}p\r)\f{\Delta u_n}{2B^{(1-(\f{\Delta}p))/2}}\cdot\f{u_{p-(\f{\Delta}p)}}p\pmod p.
\endalign$$
By Lemma 2.4,
$$\f{v_{p-(\f{\Delta}p)}-2B^{(1-(\f{\Delta}p))/2}}p\eq B^{(1-(\f{\Delta}p))/2}\f{B^{p-1}-1}p\pmod p.$$
So we have the desired (1.13) and (1.14).

In view of the above, we have completed the proof of Theorem 1.2. \qed

\heading{3. Proof of Theorem 1.1}\endheading

\proclaim{Lemma 3.1} Let $A,B\in\Z$. For any $k,l\in\N$ with $k\gs l$, we have
$$u_k(A,B)v_l(A,B)-u_l(A,B)v_k(A,B)=2B^lu_{k-l}(A,B)\tag3.1$$
and
$$v_k(A,B)v_l(A,B)-\Delta u_k(A,B)u_l(A,B)=2B^lv_{k-l}(A,B).\tag3.2$$
\endproclaim
\Proof. (i) Clearly (4.1) holds for $l=0$. If $k\in\Z^+$, then
$$\align &u_k(A,B)v_1(A,B)-u_1(A,B)v_k(A,B)
\\=&Au_k(A,B)-v_k(A,B)=Au_k(A,B)-(2u_{k+1}(A,B)-Au_k(A,B))
\\=&2(Au_k(A,B)-u_{k+1}(A,B))=2Bu_{k-1}(A,B).
\endalign$$

Now let $k\gs l\gs 2$ and assume that for each $j=1,2$ the identity (3.1) with $l$ replaced by $l-j$ still holds.
Then
$$\align &u_k(A,B)v_l(A,B)-u_l(A,B)v_k(A,B)
\\=&u_k(A,B)(Av_{l-1}(A,B)-Bv_{l-2}(A,B))-(Au_{l-1}(A,B)-Bu_{l-2}(A,B))v_k(A,B)
\\=&A(u_k(A,B)v_{l-1}(A,B)-u_{l-1}(A,B)v_k(A,B))
\\&-B(u_k(A,B)v_{l-2}(A,B)-u_{l-2}(A,B)v_k(A,B))
\\=&A\times 2B^{l-1}u_{k-(l-1)}(A,B)-B\times2B^{l-2}u_{k-(l-2)}(A,B)
\\=&2B^{l-1}(Au_{k-l+1}(A,B)-u_{k-l+2}(A,B))=2B^lu_{k-l}(A,B).
\endalign$$
This proves (3.1) by induction on $l$.

(ii) We prove (3.2) in another way. Let $\al$ and $\beta$ be the two roots of the equation $x^2-Ax+B=0$.
Then $\al\beta=B$ and $\Delta=(\al-\beta)^2$. Hence
$$\align &v_k(A,B)v_l(A,B)-\Delta u_k(A,B)u_l(A,B)
\\=&(\al^k+\beta^k)(\al^l+\beta^l)-(\al^k-\beta^k)(\al^l-\beta^l)
\\=&2(\al^k\beta^l+\al^l\beta^k)=2(\al\beta)^l(\al^{k-l}+\beta^{k-l})=2B^lv_{k-l}(A,B).
\endalign$$
This concludes the proof.  \qed

\proclaim{Lemma 3.2} For any $m\in\Z$ and $n\in\Z^+$, we have
$$\align\sum_{r=0}^{n-1}\bi{2n}ru_{n-r}(m-2,1)=&\sum_{k=0}^{n-1}\bi{2k}km^{n-1-k},\tag3.3
\\\sum_{r=0}^{n-1}\bi{2n}rv_{n-r}(m-2,1)=&m^n-\bi{2n}n,\tag3.4
\\\sum_{r=0}^{n-1}\bi{2n-1}ru_{n-r}(m-2,1)=&\f12\sum_{k=0}^{n-1}\bi{2k}km^{n-1-k}+\f{m^{n-1}}2,\tag3.5
\\\sum_{r=0}^{n-1}\bi{2n-1}rv_{n-r}(m-2,1)=&\f{m-4}2\sum_{k=0}^{n-1}\bi{2k}km^{n-1-k}+\f{m^n}2\tag3.6
\endalign$$
\endproclaim
\Proof. [ST11, (2.1)] with $d=0$ yields (3.3). Also, [ST11, (2.1)] with $d=1$ gives
$$\sum_{r=0}^n\bi{2n}ru_{n+1-r}=\sum_{k=0}^{n-1}\bi{2k}{k+1}m^{n-1-k}+m^n.\tag3.7$$
Hence
$$\align \sum_{r=0}^{n-1}\bi{2n}rv_{n-r}=&\sum_{r=0}^{n-1}\bi{2n}r(2u_{n+1-r}-(m-2)u_{n-r})
\\=&2\(\sum_{k=0}^{n-1}\bi{2k}{k+1}m^{n-1-k}+m^n-\bi{2n}n\)
\\&-(m-2)\sum_{k=0}^{n-1}\bi{2k}km^{n-1-k}
\\=&2m^n-2\bi{2n}n+\sum_{k=0}^{n-1}\l(2\bi{2k}{k+1}+(2-m)\bi{2k}k\r)m^{n-1-k}
\endalign$$
and thus (3.4) follows since
$$\align&\sum_{k=0}^{n-1}\l(2\bi{2k+1}{k+1}-m\bi{2k}k\r)m^{n-1-k}
\\=&\sum_{k=0}^{n-1}\bi{2(k+1)}{k+1}m^{n-(k+1)}-\sum_{k=0}^{n-1}\bi{2k}km^{n-k}=\bi{2n}n-m^n.
\endalign$$

Substituting $n-1$ for $n$ in (3.3) and (3.7), we get
$$\sum_{r=0}^{n-1}\bi{2(n-1)}ru_{n-1-r}(m-2,1)=\sum_{0\ls k<n-1}\bi{2k}km^{n-2-k}$$
and
$$\sum_{s=0}^{n-1}\bi{2(n-1)}su_{(n-1)+1-s}(m-2,1)=\sum_{0\ls k<n-1}\bi{2k}{k+1}m^{n-2-k}+m^{n-1}.$$
Adding the last two identities and noting that
$$\bi{2k}k+\bi{2k}{k+1}=\bi{2k+1}k=\f12\bi{2(k+1)}{k+1}\ \ \t{for all}\ k\in\N,$$
we obtain
$$\align&\sum_{0\ls k<n-1}\f12\bi{2(k+1)}{k+1}m^{n-2-k}+m^{n-1}
\\=&\sum_{r=0}^{n-1}\bi{2(n-1)}ru_{n-1-r}(m-2,1)
\\&+\sum_{0\ls r<n-1}\bi{2(n-1)}{r+1}u_{n-1-r}(m-2,1)+u_n(m-2,1)
\\=&\sum_{0\ls r<n-1}\bi{2n-1}{r+1}u_{n-1-r}(m-2,1)+u_n(m-2,1)
\endalign$$
and hence (3.5) holds.
(3.3) minus (3.5) gives
$$\sum_{0<r<n}\bi{2n-1}{r-1}u_{n-r}(m-2,1)=\f12\sum_{k=0}^{n-1}\bi{2k}km^{n-1-k}-\f{m^{n-1}}2.\tag3.8$$
By induction,
$$v_k(m-2,1)=(m-2)u_k(m-2,1)-2u_{k-1}(m-2,1)\quad\t{for all}\ k\in\Z^+.$$
Therefore
$$\align&\sum_{r=0}^{n-1}\bi{2n-1}rv_{n-r}(m-2,1)
\\=&(m-2)\sum_{r=0}^{n-1}\bi{2n-1}ru_{n-r}(m-2,1)-2\sum_{r=0}^{n-1}\bi{2n-1}{(r+1)-1}u_{n-(r+1)}(m-2,1)
\\=&\f{m-2}2\sum_{k=0}^{n-1}\bi{2k}km^{n-1-k}+\f{m-2}2m^{n-1}-\sum_{k=0}^{n-1}\bi{2k}km^{n-1-k}+m^{n-1}
\\=&\f{m-4}2\sum_{k=0}^{n-1}\bi{2k}km^{n-1-k}+\f{m^n}2
\endalign$$
with the helps of (3.5) and (3.8). This proves (3.6). \qed

\proclaim{Lemma 3.3} Let $m\in\Z\sm\{0\}$ and $\Delta=m(m-4)$. And let $p$ be an odd prime with $p\nmid \Delta$. Then
$$\sum_{k=1}^{p-1}\bi pku_k(m-2,1)\eq\l(\f{\Delta}p\r)\f{m^{p-1}-1}2+\f{4-m}4u_{p-(\f{\Delta}p)}(m-2,1)\pmod{p^2}\tag3.9$$
and
$$\sum_{k=1}^{p-1}\bi pkv_k(m-2,1)\eq\f m2(m^{p-1}-1)-\f{\Delta}4\l(\f{\Delta}p\r)u_{p-(\f{\Delta}p)}(m-2,1)
\pmod{p^2}\pmod {p^2}.\tag3.10$$
\endproclaim
\Proof. Let $\al$ and $\beta$ be the two roots of the equation $x^2-(m-2)x+1=0$. Then
$$\align &2+v_p(m-2,1)+\sum_{k=1}^{p-1}\bi pkv_k(m-2,1)
\\=&\sum_{k=0}^{p}\bi pk(\al^k+\beta^k)=(1+\al)^p+(1+\beta)^p=v_p(m,m).
\endalign$$
since
$$(1+\al)+(1+\beta)=2+m-2=m$$
 and $$(1+\al)(1+\beta)=1+(\al+\beta)+\al\beta=1+(m-2)+1=m.$$
Similarly,
$$\align &u_p(m-2,1)+\sum_{k=1}^{p-1}\bi pku_k(m-2,1)
\\=&\sum_{k=0}^p\bi pk\f{\al^k-\beta^k}{\al-\beta}
=\f{(1+\al)^p-(1+\beta)^p}{(1+\al)-(1+\beta)}=u_p(m,m).
\endalign$$

In view of [S10, Lemma 2.4],
$$2u_p(m,m)-\l(\f{\Delta}p\r)m^{p-1}\eq u_p(m-2,1)+u_{p-(\f{\Delta}p)}(m-2,1)\pmod{p^2}.\tag3.11$$
By [S12a, (3.6)] we have
$$u_p(m-2,1)-\l(\f{\Delta}p\r)\eq\f{m-2}2u_{p-(\f{\Delta}p)}(m-2,1)\pmod{p^2}.\tag3.12$$
Therefore,
$$\align 2\sum_{k=1}^{p-1}\bi pku_k(m-2,1)=&2(u_p(m,m)-u_p(m-2,1))
\\\eq&\l(\f{\Delta}p\r)m^{p-1}-u_p(m-2,1)+u_{p-(\f{\Delta}p)}(m-2,1)
\\\eq&\l(\f{\Delta}p\r)(m^{p-1}-1)+\l(\f{2-m}2+1\r)u_{p-(\f{\Delta}p)}(m-2,1)\pmod{p^2}.
\endalign$$
This proves (3.9).

By the paragraph following [S10, (2.10)],
$$\l(\f{\Delta}p\r)\f{v_p(m,m)}m\eq\l(\f{\Delta}p\r)m^{p-1}-(u_p(m,m)-u_p(m-2,1))\pmod{p^2}.\tag3.13$$
So we have
$$\align\l(\f{\Delta}p\r)\f{v_p(m,m)}m\eq&\l(\f{\Delta}p\r)\l(m^{p-1}-\f{m^{p-1}-1}2\r)-\f{4-m}4u_{p-(\f{\Delta}p)}(m-2,1)
\\=&\l(\f{\Delta}p\r)\f{m^{p-1}+1}2+\f{m-4}4u_{p-(\f{\Delta}p)}(m-2,1)\pmod{p^2}.
\endalign$$
As
$$\align v_p(m-2,1)=&2u_{p+1}(m-2,1)-(m-2)u_p(m-2,1)
\\=&(m-2)u_p(m-2,1)-2u_{p-1}(m-2,1),
\endalign$$
we have
$$\align\l(\f{\Delta}p\r)v_p(m-2,1)=&(m-2)u_p(m-2,1)-2u_{p-(\f{\Delta}p)}(m-2,1)
\\\eq&(m-2)\l(\l(\f{\Delta}p\r)+\f{m-2}2u_{p-(\f{\Delta}p)}(m-2,1)\r)-2u_{p-(\f{\Delta}p)}(m-2,1)
\\=&(m-2)\l(\f{\Delta}p\r)+\f{\Delta}2u_{p-(\f{\Delta}p)}(m-2,1)\pmod{p^2}.
\endalign$$
Combining this with (3.13), we finally get
$$\align\sum_{k=1}^{p-1}\bi pkv_k(m-2,1)=&v_p(m,m)-v_p(m-2,1)-2
\\\eq&\f m2(m^{p-1}+1)+\l(\f{\Delta}p\r)\f{\Delta}4u_{p-(\f{\Delta}p)}(m-2,1)
\\&-(m-2)-\l(\f{\Delta}p\r)\f{\Delta}2u_{p-(\f{\Delta}p)}(m-2,1)-2
\\=&\f m2(m^{p-1}-1)-\f{\Delta}4\l(\f{\Delta}p\r)u_{p-(\f{\Delta}p)}(m-2,1)
\pmod{p^2}.
\endalign$$
This proves (3.10).

In view of the above, we have completed the proof. \qed

\proclaim{Lemma 3.4} Let $m\in\Z$, and let $p$ be an odd prime $p$ not dividing $\Delta=m(m-4)$.
Then, for any $n\in\Z^+$ we have
$$\aligned&\sum_{r=0}^{n-1}\bi{2n}ru_{p(n-r)}(m-2,1)
\\\eq&\l(\l(\f{\Delta}p\r)+\f{m-4}2nu_{p-(\f{\Delta}p)}(m-2,1)\r)\sum_{k=0}^{n-1}\bi{2k}km^{n-1-k}
\\&+n\bi{2n-1}{n-1}u_{p-(\f{\Delta}p)}(m-2,1)\pmod{p^{2+\ord_p(n)}}.
\endaligned\tag3.14$$
\endproclaim
\Proof. For simplicity, we write $u_k=u_k(m-2,1)$ and $v_k=v_k(m-2,1)$ for all $k\in\N$. For each $r=0,\ldots,n$,
by Theorem 1.2 we have
$$u_{pr}\eq\l(\f{\Delta}p\r)u_r+\f{pr}2v_r\f{u_{p-(\f{\Delta}p)}}p\pmod{p^{2+\ord_p(r)}}\tag3.15$$
and
$$v_{pr}\eq v_r+\f{pr}2\Delta u_r\l(\f{\Delta}p\r)\f{u_{p-(\f{\Delta}p)}}p\pmod{p^{2+\ord_p(r)}}.\tag3.16$$
Since $r\bi{2n}r=2n\bi{2n-1}{r-1}$ for all $r\in\Z^+$, by (3.15) and (3.16) we get
$$\align &\sum_{r=0}^{n-1}\bi{2n}ru_{pr}
\\\eq&\sum_{r=0}^{n-1}\bi{2n}r\l(\l(\f{\Delta}p\r)u_r+\f{r}2v_ru_{p-(\f{\Delta}p)}\r)
\\=&\l(\f{\Delta}p\r)\sum_{r=0}^{n-1}\bi{2n}ru_r+nu_{p-(\f{\Delta}p)}\sum_{0<r<n}\bi{2n-1}{r-1}v_r\pmod{p^{2+\ord_p(n)}}
\endalign$$
and
$$\align &\sum_{r=0}^{n-1}\bi{2n}rv_{pr}
\\\eq&\sum_{r=0}^{n-1}\bi{2n}r\l(v_r+\f{r}2\Delta u_r\l(\f{\Delta}p\r)u_{p-(\f{\Delta}p)}\r)
\\=&\sum_{r=0}^{n-1}\bi{2n}rv_r+n\Delta \l(\f{\Delta}p\r)u_{p-(\f{\Delta}p)}\sum_{0<r<n}\bi{2n-1}{r-1}u_r\pmod{p^{2+\ord_p(n)}}.
\endalign$$
Therefore
$$\align&\sum_{r=0}^{n-1}\bi{2n}r(u_{pn}v_{pr}-v_{pn}u_{pr})
\\\eq&u_{pn}\sum_{r=0}^{n-1}\bi{2n}rv_r-v_{pn}\l(\f{\Delta}p\r)\sum_{r=0}^{n-1}\bi{2n}ru_r
\\&+u_{pn}\times n\Delta\l(\f{\Delta}p\r)u_{p-(\f{\Delta}p)}\sum_{0<r<n}\bi{2n-1}{r-1}u_r
\\&-v_{pn}\times nu_{p-(\f{\Delta}p)}\sum_{0<r<n}\bi{2n-1}{r-1}v_r\pmod{p^{2+\ord_p(n)}}.
\endalign$$
In view of (3.15) and (3.16) with $r=n$, from the above we have
$$\align&\sum_{r=0}^{n-1}\bi{2n}r(u_{pn}v_{pr}-v_{pn}u_{pr})
\\\eq&\l(\l(\f{\Delta}p\r)u_n+\f n2v_nu_{p-(\f{\Delta}p)}\r)\sum_{r=0}^{n-1}\bi{2n}rv_r
\\&-\l(v_n+\f n2\Delta u_n\l(\f{\Delta}p\r)u_{p-(\f{\Delta}p)}\r)\l(\f{\Delta}p\r)\sum_{r=0}^{n-1}\bi{2n}ru_r
\\&+nu_n\Delta u_{p-(\f{\Delta}p)}\sum_{0<r<n}\bi{2n-1}{r-1}u_r
\\&-nv_n u_{p-(\f{\Delta}p)}\sum_{0<r<n}\bi{2n-1}{r-1}v_r
\pmod{p^{2+\ord_p(n)}}\endalign$$
and hence
$$\align&\sum_{r=0}^{n-1}\bi{2n}r(u_{pn}v_{pr}-v_{pn}u_{pr})
\\\eq&\l(\f{\Delta}p\r)\sum_{r=0}^{n-1}\bi{2n}r(u_nv_r-v_nu_r)
\\&+\f n2u_{p-(\f{\Delta}p)}\sum_{r=0}^{n-1}\bi{2n}r(v_nv_r-\Delta u_nu_r)
\\&-nu_{p-(\f{\Delta}p)}\sum_{0<r<n}\bi{2n-1}{r-1}(v_nv_r-\Delta u_nu_r)
\pmod{p^{2+\ord_p(n)}}.
\endalign$$
Thus, with the help of Lemma 3.1 we obtain
$$\align \sum_{r=0}^{n-1}\bi{2n}ru_{pn-pr}
\eq&\l(\f{\Delta}p\r)\sum_{r=0}^{n-1}\bi{2n}ru_{n-r}
+\f n2u_{p-(\f{\Delta}p)}\sum_{r=0}^{n-1}\bi{2n}rv_{n-r}
\\&-nu_{p-(\f{\Delta}p)}\sum_{r=0}^{n-1}\l(\bi{2n}r-\bi{2n-1}r\r)v_{n-r}
\\\eq&nu_{p-(\f{\Delta}p)}\(\sum_{r=0}^{n-1}\bi{2n-1}rv_{n-r}-\f12\sum_{r=0}^{n-1}\bi{2n}rv_{n-r}\)
\\&+\l(\f{\Delta}p\r)\sum_{r=0}^{n-1}\bi{2n}ru_{n-r}\pmod{p^{2+\ord_p(n)}}.
\endalign$$
Combining this with Lemma 3.2, we immediately get (3.14). \qed

\proclaim{Lemma 3.5} For any $m\in\Z\sm\{0\}$ and $n\in\Z^+$, we have
$$\sum_{k=0}^{n-1}\f{\bi{2k}{k+1}}{m^k}=\f{\bi{2n-1}{n-1}}{m^{n-1}}-\f m2+\f{m-2}2\sum_{k=0}^{n-1}\f{\bi{2k}k}{m^k}.\tag3.17$$
\endproclaim
\Proof. Observe that
$$\align \sum_{k=0}^{n-1}\f{\bi{2k}k+\bi{2k}{k+1}}{m^k}=&\sum_{k=0}^{n-1}\f{\bi{2k+1}k}{m^k}=\f{\bi{2n-1}{n-1}}{m^{n-1}}+\f m2\sum_{0\ls k<n-1}\f{\bi{2k+2}{k+1}}{m^{k+1}}
\\=&\f{\bi{2n-1}{n-1}}{m^{n-1}}+\f m2\(\sum_{r=0}^{n-1}\f{\bi{2r}r}{m^r}-1\).
\endalign$$
So (3.17) follows. \qed

\medskip
\noindent{\it Proof of Theorem 1.1}. We first handle the case $\Delta=m(m-4)\eq0\pmod p$.
By [S11a, Theorem 1.1],
$$\f1{pn}\sum_{k=0}^{pn-1}\f{\bi{2k}k}{m^k}\eq\f{\bi{2pn-1}{pn-1}}{4^{pn-1}}+\da_{p,3}\f{m-4}3\bi{2n/p^{\ord_p(n)}-1}{n/p^{\ord_p(n)}-1}\pmod p.$$
By Lucas' theorem.
$$\bi{2pn-1}{pn-1}=\f12\bi{2pn}{pn}\eq\f12\bi{2n}n\eq\f12\bi{2n/p^{\ord_p(n)}}{n/p^{\ord_p(n)}}=\bi{2n/p^{\ord_p(n)}-1}{n/p^{\ord_p(n)}-1}\pmod p.$$
Since $m\eq4\pmod p$, we have $m\eq1\pmod p$ if $p=3$. Therefore
$$\align\f1{pn}\sum_{k=0}^{pn-1}\f{\bi{2k}k}{m^k}
\eq&\f{\bi{2n-1}{n-1}}{4^{pn-1}}+\da_{p,3}\f{m-4}3\bi{2n-1}{n-1}
\\\eq&\f{\bi{2n-1}{n-1}}{4^{n-1}}+\da_{p,3}\f{m(m-4)}{3m^{n-1}}\bi{2n-1}{n-1}
\\\eq&\f{\bi{2n-1}{n-1}}{m^{n-1}}\l(1+\da_{p,3}\f{\Delta}3\r)\pmod p.
\endalign$$
By Lemma 2.2,
$$\f{u_p(m-2,1)}p\eq\l(\f{m-2}2\r)^{p-1}+\da_{p,3}\f{\Delta}3\eq1+\da_{p,3}\f{\Delta}3\pmod p.$$
So
$$\f1n\sum_{k=0}^{pn-1}\f{\bi{2k}k}{m^k}\eq\f{\bi{2n-1}{n-1}}{m^{n-1}}u_p(m-2,1)\pmod{p^2}$$
as desired.

Below we assume that $p\nmid \Delta$.
By [ST11, Lemma 2.1], for any $r=1,\ldots,n$ we have
$$\bi{2pn}{pr}\big/\bi{2n}r\eq1\pmod{p^{2+\ord_p(r)}}$$
and hence
$$\bi{2pn}{pr}\eq\bi{2n}r=\f{2n}r\bi{2n-1}{r-1}\pmod{p^{2+\ord_p(n)}}.\tag3.18$$
By Lucas' theorem, for any $r\in\N$ and $k\in\{1,\ldots,p\}$ we have
$$\bi{p(2n-1)+p-1}{pr+k-1}\eq \bi{2n-1}r\bi{p-1}{k-1}\pmod p.$$
So, in view of Lemmas 3.1, 3.2 and 3.4, we have
$$\align&\sum_{k=0}^{pn-1}\bi{2k}km^{pn-1-k}
\\=&\sum_{s=0}^{pn-1}\bi{2pn}su_{pn-s}=\sum_{r=0}^{n-1}\sum_{k=0}^{p-1}\bi{2pn}{pr+k}u_{pn-(pr+k)}
\\=&\sum_{r=0}^{n-1}\bi{2pn}{pr}u_{p(n-r)}
+\sum_{r=0}^{n-1}\sum_{k=1}^{p-1}\f{2pn}{pr+k}\bi{p(2n-1)+p-1}{pr+k-1}u_{p(n-r)-k}
\\\eq&\sum_{r=0}^{n-1}\bi{2n}ru_{p(n-r)}+\sum_{r=0}^{n-1}\sum_{k=1}^{p-1}\f{2pn}k\bi{2n-1}r\bi{p-1}{k-1}u_{p(n-r)-k}
\\\eq&\l(\l(\f{\Delta}p\r)+\f{m-4}2nu_{p-(\f{\Delta}p)}\r)\sum_{k=0}^{n-1}\bi{2k}km^{n-1-k}
+n\bi{2n-1}{n-1}u_{p-(\f{\Delta}p)}
\\&+n\sum_{r=0}^{n-1}\bi{2n-1}r\sum_{k=1}^{p-1}\bi pk(u_{p(n-r)}v_k-u_kv_{p(n-r)})\pmod{p^{2+\ord_p(n)}}.
\endalign$$
and hence
$$\align\sum_{k=0}^{pn-1}\f{\bi{2k}k}{m^k}
\eq&\l(\l(\f{\Delta}p\r)m^{(1-p)n}+\f{m-4}2nu_{p-(\f{\Delta}p)}\r)\sum_{r=0}^{n-1}\f{\bi{2r}r}{m^r}
+\f{n}{m^{n-1}}\bi{2n-1}{n-1}u_{p-(\f{\Delta}p)}
\\&+\f{n}{m^{n-1}}\sum_{r=0}^{n-1}\bi{2n-1}r\sum_{k=1}^{p-1}\bi pk(u_{p(n-r)}v_k-u_kv_{p(n-r)})
\pmod{p^{2+\ord_p(n)}}.
\endalign$$
By Lemma 2.5,
$$\f1{m^{(p-1)n}}\eq\f1{1+n(m^{p-1}-1)}\eq1-n(m^{p-1}-1)\pmod{p^{2+\ord_p(n)}}.\tag3.19$$
Therefore
$$\aligned&\f1n\(\sum_{k=0}^{pn-1}\f{\bi{2k}k}{m^k}-\l(\f{\Delta}p\r)\sum_{r=0}^{n-1}\f{\bi{2r}r}{m^r}\)-\f{\bi{2n-1}{n-1}}{m^{n-1}}u_{p-(\f{\Delta}p)}
\\\eq&\l(\l(\f{\Delta} p\r)(1-m^{p-1})+\f{m-4}2u_{p-(\f{\Delta}p)}\r)\sum_{r=0}^{n-1}\f{\bi{2r}r}{m^r}
\\&+\f{1}{m^{n-1}}\sum_{r=0}^{n-1}\bi{2n-1}r\sum_{k=1}^{p-1}\bi pk(u_{p(n-r)}v_k-u_kv_{p(n-r)})
\pmod{p^2}.
\endaligned\tag3.20$$

Note that $p\mid \bi pk$ for all $k=1,\ldots,p-1$. In light of Theorem 1.2 and (3.5)-(3.6),
$$\align&\sum_{r=0}^{n-1}\bi{2n-1}r\sum_{k=1}^{p-1}\bi pk(u_{p(n-r)}v_k-u_kv_{p(n-r)})
\\\eq&\sum_{r=0}^{n-1}\bi{2n-1}r\sum_{k=1}^{p-1}\bi pk\l(\l(\f{\Delta}p\r)u_{n-r}v_k-u_kv_{n-r}\r)
\\=&\l(\f{\Delta}p\r)\sum_{k=1}^{p-1}\bi pkv_k\sum_{r=0}^{n-1}\bi{2n-1}ru_{n-r}-\sum_{k=1}^{p-1}\bi pku_k\sum_{r=0}^{n-1}\bi{2n-1}rv_{n-r}
\\=&\l(\f{\Delta}p\r)\sum_{k=1}^{p-1}\bi pkv_k\(\f12\sum_{r=0}^{n-1}\bi{2r}rm^{n-1-r}+\f{m^{n-1}}2\)
\\&-\sum_{k=1}^{p-1}\bi pku_k\(\f{m-4}2\sum_{r=0}^{n-1}\bi{2r}rm^{n-1-r}+\f{m^n}2\)\pmod{p^2}.
\endalign$$
Combining this with (3.20) we have reduced (1.5) to the congruence
$$\align&\l(\l(\f{\Delta}p\r)(m^{p-1}-1)+\f{4-m}2u_{p-(\f{\Delta}p)}\r)\sum_{r=0}^{n-1}\f{\bi{2r}r}{m^r}
\\\eq&\l(\f{\Delta}p\r)\sum_{k=1}^{p-1}\bi pk\f{v_k}2\(\sum_{r=0}^{n-1}\f{\bi{2r}r}{m^r}+1\)
\\&-\sum_{k=1}^{p-1}\bi pk\f{u_k}2\((m-4)\sum_{r=0}^{n-1}\f{\bi{2r}r}{m^r}+m\)\pmod{p^2}.
\endalign$$
This indeed holds since
$$\l(\f{\Delta}p\r)\sum_{k=1}^{p-1}\bi pkv_k\eq m\sum_{k=1}^{p-1}\bi pku_k\pmod{p^2}$$
and
$$2\sum_{k=1}^{p-1}\bi pku_k\eq \l(\f{\Delta}p\r)(m^{p-1}-1)+\f{4-m}2u_{p-(\f{\Delta}p)}\pmod{p^2}$$
by Lemma 3.3. So we finally obtain (1.5).

Next we deduce (1.6) and (1.7) via (1.5). By Lemma 3.5,
$$\sum_{k=0}^{pn-1}\f{\bi{2k}{k+1}}{m^k}=\f{\bi{2pn-1}{pn-1}}{m^{pn-1}}-\f m2+\f{m-2}2\sum_{k=0}^{pn-1}\f{\bi{2k}k}{m^k}.\tag3.21$$
Since
$$\bi{2pn-1}{pn-1}=\f12\bi{2pn}{pn}\eq\f12\bi{2n}n=\bi{2n-1}{n-1}\pmod{p^{2+\ord_p(n)}}$$
by (3.18), from (3.19), (3.21) and (1.5) we get
$$\align&\sum_{k=0}^{pn-1}\f{\bi{2k}{k+1}}{m^k}-\f{\bi{2n-1}{n-1}}{m^{n-1}}(1-n(m^{p-1}-1))+\f m2
\\\eq&\f{m-2}2\l(\f{\Delta}p\r)\sum_{r=0}^{n-1}\f{\bi{2r}r}{m^r}+\f{m-2}2\cdot\f{n}{m^{n-1}}\bi{2n-1}{n-1}u_{p-(\f{\Delta}p)}\pmod{p^{2+\ord_p(n)}}.
\endalign$$
Combining this with (3.17) we immediately obtain (1.6).
As $C_k=\bi{2k}k-\bi{2k}{k+1}$ for $k\in\N$, (1.7) follows from (1.5) and (1.6).

The proof of Theorem 1.1 is now complete. \qed

\heading{4. Proof of Theorem 1.3}\endheading

\proclaim{Lemma 4.1} Let $p>3$ be a prime. Then, for any integers $n\gs k\gs0$, we have
$$\f{\bi{pn}{pk}}{\bi nk}\in 1+p^3nk(n-k)\Z_p.\tag4.1$$
\endproclaim

\Remark\ 4.1. This is a useful known result, see, e.g., [RZ].

\proclaim{Lemma 4.2} Let $p$ be a prime, and let $j\in\N$ and $k\in\{0,\ldots,p-1\}$. Then
$$\bi{2jp+2k}{jp+k}\eq\bi{2j}j\bi{2k}k\pmod p.\tag4.2$$
\endproclaim
\Proof. If $2k<p$ then (4.2) follows from Lucas' congruence. If $2k\gs p$, then $p\mid\bi{2k}k$, and
by the Lucas congruence we have
$$\align\bi{2jp+2k}{jp+k}=&\bi{(2j+1)p+(2k-p)}{jp+k}
\\\eq&\bi{2j+1}j\bi{2k-p}k=0\eq\bi{2j}j\bi{2k}k\pmod p.
\endalign$$
This concludes the proof. \qed

\medskip
\noindent{\it Proof of Theorem 1.3}. Since
 $$\bi nk^2nk(n-k)=n^3\bi{n-1}{k-1}\bi{n-1}{n-k-1}$$
 for any positive integer $k<n$, by Lemma 4.1 we have
 $$\bi{pn}{pk}^2-\bi nk^2\in p^3n^{3}\Z_p\quad\t{for all}\ k=0,\ldots,n.$$
 Thus
$$\align \sum_{k=0}^n\bi{pn}{pk}^2\bi{2pk}{pk}(-1)^{pk}
 \eq\sum_{k=0}^n\bi nk^2\bi{2pk}{pk}(-1)^k\pmod{p^{3+3\ord_p(n)}}.
 \endalign$$
For each $k=1,2,3,\ldots$, clearly
$$\bi{2pk}{pk}-\bi{2k}k\in p^3k^3\Z_p$$
by Lemma 4.1, and
$$\bi nk^2k^2=n^2\bi{n-1}{k-1}^2\eq0\pmod{n^2}.$$
Therefore,
$$\align\sum^{pn}\Sb k=0\\p\mid k\endSb\bi{pn}k^2\bi{2k}k(-1)^k
=&\sum_{k=0}^n\bi{pn}{pk}^2\bi{2pk}{pk}(-1)^{pk}
\\ \eq&\sum_{k=0}^n\bi nk^2\bi{2k}{k}(-1)^k=g_n(-1)\pmod{p^{3+2\,\ord_p(n)}}.
\endalign$$

In view of Lemma 4.2,
$$\align&\sum^n\Sb k=0\\p\nmid k\endSb\bi{pn}k^2\bi{2k}k(-1)^k
\\=&\sum_{r=0}^{n-1}\sum_{k=1}^{p-1}\f{(pn)^2}{(rp+k)^2}\bi{(n-1)p+p-1}{rp+k-1}^2\bi{2rp+2k}{rp+k}(-1)^{rp+k}
\\\eq&\sum_{r=0}^{n-1}\sum_{k=1}^{p-1}\f{(pn)^2}{k^2}\bi{n-1}r^2\bi{p-1}{k-1}^2\bi{2r}r\bi{2k}k(-1)^{r+k}
\\=&(pn)^2\sum_{r=0}^{n-1}\bi{n-1}r^2\bi{2r}r(-1)^r\sum_{k=1}^{p-1}\f{(-1)^{(k-1)2+k}}{k^2}\bi{2k}k
\pmod{p^{3+2\,\ord_p(n)}}.
\endalign$$
Note that
$$\sum_{k=1}^{p-1}\f{\bi{2k}k}k\eq0\pmod {p}$$
by [PS] or [ST10], and that 
$$\sum_{k=1}^{p-1}\f{(-1)^k}{k^2}\bi{2k}k\eq0\pmod p$$
by [T].

Combining the above arguments, we get
$$\align g_{pn}(-1)=&\sum_{k=0}^{pn}\bi{pn}k^2\bi{2k}k(-1)^k
\\\eq&\sum_{k=0}^n\bi nk^2\bi{2k}k(-1)^k=g_n(-1)\pmod{p^{3+2\ord_p(n)}}
\endalign$$
So $(1.15)$ holds. This concludes our proof of Theorem 1.3. \qed

\heading{5. Some conjectures}\endheading

\proclaim{Conjecture 5.1} Let $p$ be an odd prime. For any integer $m\not\eq0\pmod p$ and positive integer $n$, we have
$$\f1{n\bi{2n-1}{n-1}}\(\sum_{k=0}^{pn-1}\f{\bi{2k}k}{m^k}-\l(\f{\Delta}p\r)\sum_{r=0}^{n-1}\f{\bi{2r}r}{m^r}\)
\eq\f{u_{p-(\f{\Delta}p)}(m-2,1)}{m^{n-1}}\pmod{p^2},\tag5.1$$
where $\Delta=m(m-4)$.
\endproclaim
\Remark\ 5.1. (5.1) is stronger than (1.5).

\proclaim{Conjecture 5.2} {\rm (i)} Let $p$ be an odd prime. For any $n\in\Z^+$ we have
$$\f{\sum_{k=0}^{pn-1}\bi{2k}k-(\f p3)\sum_{r=0}^{n-1}\bi{2r}{r}}{n^2\bi{2n-1}{n-1}}\eq\sum_{k=0}^{p-1}\bi{2k}k-\l(\f p3\r)\pmod{p^4}.\tag5.2$$

{\rm (ii)} Let $p>3$ be a prime, and let $m\in\{2,3\}$ and $\Delta=m(m-4)$. Then there is a $p$-adic integer $c_p^{(m)}$ only depending on $p$ and $m$ such that for any $n\in\Z^+$ we have
$$\f{m^{n-1}}{n^2\bi{2n-1}{n-1}}\(\sum_{k=0}^{pn-1}\f{\bi{2k}k}{m^k}-\l(\f{\Delta}p\r)\sum_{r=0}^{n-1}\f{\bi{2r}r}{m^r}\)
\eq\sum_{k=0}^{p-1}\f{\bi{2k}k}{m^k}-\l(\f{\Delta}p\r)+p^3c_p^{(m)}(n-1)\pmod{p^4}.\tag5.3$$
\endproclaim
\Remark\ 5.2. In 1992 N. Strauss, J. Shallit and D. Zagier [SSZ] proved that
$$\f{\sum_{k=0}^{n-1}\bi{2k}k}{n^2\bi{2n}n}\eq-1\pmod3\quad\t{for any}\ n\in\Z^+.$$
In 2011 the author [S11b] showed that
$$\sum_{k=0}^{p-1}\f{\bi{2k}k}{2^k}\eq\l(\f{-1}p\r)-p^2E_{p-3}\pmod {p^3}\quad \ \t{for any odd prime}\ p,$$
where $E_0,E_1,E_2,\ldots$ are the Euler numbers defined by
$$\f{2}{e^x+e^{-x}}=\sum_{n=0}^\infty E_n\f{x^n}{n!}\ \ (|x|<\pi).$$

\proclaim{Conjecture 5.3} {\rm (i)} Let $p$ be an odd prime. For any integer $m\not\eq0\pmod p$ and positive integer $n$, we have
$$\f1{pn}\(\sum_{k=0}^{pn-1}\bi{pn-1}k\f{\bi{2k}k}{(-m)^k}-\l(\f{m(m-4)}p\r)\sum_{r=0}^{n-1}\bi{n-1}r\f{\bi{2r}r}{(-m)^r}\)\in\Z_p.\tag5.4$$

{\rm (ii)} For any prime $p\not=3$ and $n\in\Z^+$, we have
$$\f1{p^2n^2}\(\sum_{k=0}^{pn-1}\bi{pn-1}k\f{\bi{2k}k}{(-3)^k}-\l(\f p3\r)\sum_{r=0}^{n-1}\bi{n-1}r\f{\bi{2r}r}{(-3)^r}\)\in\Z_p.\tag5.5$$
\endproclaim
\Remark\ 5.3. The author [S12a] determined $\sum_{k=0}^{p-1}\bi{p-1}k\bi{2k}k/(-m)^k$ modulo $p^2$ for any odd prime $p$ and integer $m\not\eq0\pmod p$.

\proclaim{Conjecture 5.4} Let $p>3$ be a prime and let $n\in\Z^+$. Then
$$\align\f{16^n}{n^2\bi{2n}n^2}\(\sum_{k=0}^{pn-1}\f{\bi{2k}k^2}{16^k}-\l(\f{-1}p\r)\sum_{r=0}^{n-1}\f{\bi{2r}r^2}{16^r}\)\eq&-4p^2E_{p-3}\pmod{p^3},\tag5.6
\\\f{27^n}{n^2\bi{2n}n\bi{3n}n}\(\sum_{k=0}^{pn-1}\f{\bi{2k}k\bi{3k}k}{27^k}-\l(\f{p}3\r)\sum_{r=0}^{n-1}\f{\bi{2r}r\bi{3r}r}{27^r}\)\eq&-\f 32p^2B_{p-2}\l(\f13\r)\pmod{p^3},\tag5.7
\\\f{64^n}{n^2\bi{4n}{2n}\bi{2n}n}\(\sum_{k=0}^{pn-1}\f{\bi{4k}{2k}\bi{2k}k}{64^k}-\l(\f{-2}p\r)\sum_{r=0}^{n-1}\f{\bi{4r}{2r}\bi{2r}r}{64^r}\)\eq&-p^2E_{p-3}\l(\f14\r)\pmod{p^3},\tag5.8
\\\f{432^n}{n^2\bi{6n}{3n}\bi{3n}n}\(\sum_{k=0}^{pn-1}\f{\bi{6k}{3k}\bi{3k}k}{432^k}-\l(\f{-1}p\r)\sum_{r=0}^{n-1}\f{\bi{6r}{3r}\bi{3r}r}{432^r}\)\eq&-20p^2E_{p-3}\pmod{p^3},\tag5.9
\endalign$$
where $E_{p-3}(x)$ is the Euler polynomial of degree $p-3$.
\endproclaim
\Remark\ 5.4.  Let $p>3$ be a prime. Recently, J.-C. Liu [L] proved that for any $n\in\Z^+$ we have
$$\align\sum_{k=0}^{pn-1}\f{\bi{2k}k^2}{16^k}\eq&\l(\f{-1}p\r)\sum_{r=0}^{n-1}\f{\bi{2r}r^2}{16^r}\pmod{p^2},
\\\sum_{k=0}^{pn-1}\f{\bi{2k}k\bi{3k}k}{27^k}\eq&\l(\f{p}3\r)\sum_{r=0}^{n-1}\f{\bi{2r}r\bi{3r}r}{27^r}\pmod{p^2},
\\\sum_{k=0}^{pn-1}\f{\bi{4k}{2k}\bi{2k}k}{64^k}\eq&\l(\f{-2}p\r)\sum_{r=0}^{n-1}\f{\bi{4r}{2r}\bi{2r}r}{64^r}\pmod{p^2},
\\\sum_{k=0}^{pn-1}\f{\bi{6k}{3k}\bi{3k}k}{432^k}\eq&\l(\f{-1}p\r)\sum_{r=0}^{n-1}\f{\bi{6r}{3r}\bi{3r}r}{432^r}\pmod{p^2},
\endalign$$
the case $n=1$ of which is a conjecture of F. Rodriguez-Villegas [RV] first confirmed by E. Morterson [M03].
In the case $n=1$, (5.6) was established by the author [S11b], and (5.7)-(5.9) were conjectured by the author and later confirmed by Z.-H. Sun [Su].

\proclaim{Conjecture 5.5} For any odd prime $p$ and positive integer $n$, we have
$$\f1{n^3\bi{2n}n^3}\(\f1p\sum_{k=0}^{pn-1}(21k+8)\bi{2k}k^3-\sum_{r=0}^{n-1}(21r+8)\bi{2r}r^3\)
\eq0\pmod{p^3}.\tag5.10$$
\endproclaim
\Remark\ 5.5. (5.10) in the case $n=1$ was proved by the author in [S11b].
We guess that all those Ramanujan-type supercongruences should have extensions involving $n\in\Z^+$ similar to (5.10).
\medskip

\proclaim{Conjecture 5.6} For any prime $p>5$ and $n\in\Z^+$, we have
$$\f{g_{pn}(-1)-g_n(-1)}{(pn)^3}\in\Z_p.\tag5.11$$
\endproclaim
\Remark\ 5.6. This is stronger than Theorem 1.3. We note the new identity
$$g_n(-1)=\sum_{k=0}^{\lfloor n/2\rfloor}\bi n{2k}^2\bi{2k}k(-1)^{n-k},$$
which can be easily proved via the Zeilberger algorithm (cf. [PWZ, pp.\,101-119]).
\medskip

In the same spirit, we have many other conjectures similar to the above ones (see, e.g., 
Conjectures 12, 22-24, 26-32, 60-63 and 82 of [S19]).
In our opinion, almost all previous known congruences should have such extensions involving a parameter $n\in\Z^+$.

\enddocument